\documentclass[12pt]{article}

\usepackage{amsfonts}
\usepackage{amsmath}
\usepackage{amssymb}
\usepackage{amsthm}
\usepackage{graphics}

\newtheorem{thm}{Theorem}[section]

\newtheorem{lem}[thm]{Lemma}

\theoremstyle{plain}

\newtheorem{claim}[thm]{Claim}
\newtheorem{obs}[thm]{Observation}

\newtheorem{rem}[thm]{Remark}

\theoremstyle{remark}

\def\al{\alpha}

\def\be{\beta}

\def\de{{\delta}}

\def\Ga{\Gamma}
\def\si{\sigma}
\def\de{{\delta}}

\def\ze{\zeta}
\def\om{\omega}
\def\Om{\Omega}

\def\ra{\rightarrow}
\def\nra{\not \ra}

\def \Z {\mathbb Z}

\def\rar{\leftrightarrow}
\def\nrar{\hspace{.03in}\not\hspace{-.02in}\leftrightarrow}
\newcommand{\sub}[0]{\subseteq}
\newcommand{\sm}[0]{\setminus}
\newcommand{\0}[0]{\emptyset}
\def \mn {\medskip\noindent}
\def \bn {\bigskip\noindent}
\newcommand{\E}[0]{{\sf E}}
\newcommand{\beq}[1]{\begin{equation}\label{#1}}
\newcommand{\enq}[0]{\end{equation}}

\renewcommand{\phi}[0]{\varphi}
\newcommand{\hap}[0]{\hat{\phi}}
\newcommand{\hOm}[0]{\hat{\Om}}
\renewcommand{\dots}[0]{,\ldots ,}

\begin{document}
%%%
\renewcommand{\thefootnote}{\fnsymbol{footnote}}
%%%
\title{Some Conditional Correlation Inequalities for Percolation and
Related Processes}

\author{J. van den Berg, O. H\"{a}ggstr\"{o}m and
J. Kahn\footnotemark \\
{\small CWI, Chalmers University of Technology, and Rutgers University } \\
{\footnotesize email: J.van.den.Berg@cwi.nl; olleh@math.chalmers.se;
jkahn@math.rutgers.edu}
}
\date{}

\footnotetext{ * Supported by NSF grant DMS0200856.}

\maketitle

\begin{abstract}
Consider ordinary bond percolation on a finite or countably infinite graph.
Let $s$, $t$, $a$ and $b$ be vertices.
An earlier paper \cite{BeKa} proved the (nonintuitive)
result that, conditioned on the event that there is no open path from
$s$ to $t$,
the two events ``there is an open path from $s$ to  $a$" and
``there is an open path from $s$ to $b$" are positively correlated.
In the present paper we further investigate and generalize the
theorem of which this result
was a consequence. This leads to results saying, informally,
that, with the above conditioning, the open cluster of $s$ is
conditionally positively (self-)associated and that it
is conditionally negatively correlated with the open cluster of $t$.
%In particular, the two events
%``there is an open path from $s$ to $a$"
%and ``there is an open path from $t$ to $b$" are conditionally
%negatively correlated.
%Although this statement is intuitively correct,
%we don't know a straightforward intuitive
%proof.

We also present analogues of some of our
results for
(a) random-cluster measures, and
(b) directed percolation
and contact processes,
and observe that the latter lead to
improvements of some of the results in a paper of
Belitsky, Ferrari, Konno and
Liggett (1997).
\end{abstract}

\section{Introduction and results for ordinary percolation}
\label{Intro}
This paper is concerned with positive and negative correlation
and the stronger notion of positive association.
Recall that events $A,B$ (in some probability space) are
{\em positively} {\em correlated} if
$\Pr(AB)\geq \Pr(A)\Pr(B)$,
and negatively correlated if the reverse inequality holds.
Positive association will be defined below 
(following Theorem 1.3).

We begin in this section with results for ordinary bond percolation.
Our original motivation here (and for the present work) was 
Theorem 1.4.

We then consider extensions to the random cluster model (Section 2)
and to percolation on directed graphs, together with applications to
the contact process (Section 3).

A few words about proofs may be in order.  
The approach given for
percolation in Section 1 is similar to that of 
\cite{BeKa} (see the
proof of the present Theorem~\ref{1.1}).
This approach does not seem applicable to 
the random cluster model,
and Section 2 takes a completely different route, based
on Markov chains, to extend the results
of Section 1 to this more general setting.
We also describe, in Section~\ref{Fuzzy}, a different 
way of getting at some of the random cluster results.
This is based on a connection with the ``fuzzy Potts model," 
and is included here despite handling 
only a subset of what's
covered by the Markov chain approach, because we
think the relevance of the connection is interesting.
The results for ``directed percolation" in 
Section~\ref{DPCP} can again be obtained using either
the approach of Section~\ref{Intro} or the 
Markov chain approach of Section~\ref{RCM}.
Here we have tried to keep the discussion brief,
mainly indicating those points where the generalization
from what came before may not be entirely obvious.

\bigskip
Consider bond percolation on a (finite or countably infinite, locally finite)
graph $G = (V,E)$,
where each edge $e$ is, independently of all other edges, open with
probability $p_e$ and
closed with probability $1 - p_e$.
For $a, b \in V$ the event that there is an open path from $a$ to $b$
is denoted by
$a \rar b$, and the
complement of this event by $a \nrar b$.
For $X, Y \subset V$ we write $X \nrar Y$ for the event
$\{x \nrar y \,\, \forall x \in X, y \in Y\}$.

\medskip
In an earlier paper \cite{BeKa} we showed
that, for any vertices $s, t,a, b $,
\begin{equation}
\label{vdBK1}
\Pr(s \rar a, \,\, s \rar b \, | \, s \nrar t) \geq
\Pr(s \rar a \, \mid \, s \nrar t) \,
\Pr(s \rar b \, \mid \, s \nrar t).
\end{equation}
This was a consequence (really a special case)
of Theorem 1.2 of \cite{BeKa}, 
to which we will return below.

Here we show, among other results, a sort of complement
of \eqref{vdBK1}, {\em viz.}
\begin{equation}
\label{new1}
\Pr(s \rar a, \,\, t \rar b \, | \, s \nrar t) \leq
\Pr(s \rar a \, \mid \, s \nrar t) \,
\Pr(t \rar b \, \mid \, s \nrar t).
\end{equation}
In this section we will prove the quite intuitive \eqref{new1}
by way of a generalization of
the {\em not} very intuitive
\eqref{vdBK1}.
Before giving this generalization, we need some further
definitions and notation.

Let $s$ be a fixed vertex.
By the {\em open cluster}, $C_s$, of $s$ we mean the set of
all edges which are in open paths
starting at $s$.
As in \cite{BeKa} we define, for $X\sub V$,
the event
$$R_X := \{s \nrar X\} = \{s \nrar x \,\, \forall x \in X\}.$$

Let $\Om = \{0,1\}^E$ be the set of realizations;
elements of $\Om$ will typically be denoted $\om$.
Recall that an event $A$ is {\em increasing}
(really, {\em nondecreasing})
if $\om'\geq\om \in A$
implies
$\om' \in A$.
We also say that $A$ is {\em increasing
and determined by the open cluster of $s$} if
$\om \in A$ and $C_s(\om') \supseteq
C_s(\om)$ imply $\om' \in A$.
(Note that such an event is increasing in the sense
above.)
A simple example of such an event is $\{s \rar a\}$.

The following statement is a natural generalization of Theorem 1.2
of \cite{BeKa}.

\begin{thm}\label{1.1}
Let $A$ and $B$ be increasing events determined by the open cluster of
$s$. Then
for all $X, Y \sub V\sm\{s\}$,
\begin{equation}
\label{prop1}
\Pr(A \, R_X) \Pr(B \, R_Y)
\leq \Pr(A \, B \, R_{X \cap Y})
\Pr(R_{X \cup Y}).
\end{equation}
\end{thm}

\medskip\noindent
{\bf Remark.} Theorem 1.2 in \cite{BeKa} is the special case where
each of $A$,
$B$ is of
the form $\{s \rar w \,\, \forall w \in W\}$ for some $W \subset V$.
%For present purposes we
%need the more general version given here.
The proof of the present more general result is almost the
same and we present it
in a slightly abbreviated form, emphasizing the parts which need extra
attention because of the
generalization.  (One may say that the key idea (in
both cases) is 
generalizing from statements like (\ref{vdBK1})
to the
form (\ref{prop1}), which supports an inductive proof.)

\mn
{\em Proof.}
We give the proof for finite $G$; the infinite
case then
follows by standard limit arguments. The proof is by induction
on the number of vertices.
When $G$ has only one vertex, the result is obvious;
so we suppose, for some
$n \geq 1$, that the result holds for graphs with
at most $n$ vertices, and consider $G$ with
$n+1$ vertices.

With notation as in the theorem,
it is easy to see that there is
an event $\tilde{A} \sub A$ with the
following properties: it is increasing
and determined by $C_s$; it does not depend on
$$E_X:=\{e\in E:e\cap X\neq\0\}$$
(that is, if $\om'_e = \om_e$ for all $e\not\in E_X$,
then
$\om \in \tilde{A}$ iff $\om' \in \tilde{A}$);
and, finally, $\tilde{A} \, R_X = A \, R_X$.
A similar remark holds for $B$ and $Y$.
So we may assume that $A$ does not
depend on $E_X$ and $B$ does not depend on $E_Y$.

If $X \cap Y = \emptyset$, the r.h.s. of \eqref{prop1} is
$\Pr(A \,B)$ $\Pr(R_{X \cup Y})$, and two applications of
the FKG Inequality give the result:
\begin{equation}\label{base}
\Pr(A \, R_X) \Pr(B \, R_Y) \leq
\Pr(A)  \Pr(B)  \Pr(R_X)  \Pr(R_Y)
\leq \Pr(A \, B )
\Pr(R_{X \cup Y})
\end{equation}
(note $R_{X\cup Y}=R_XR_Y$).

Now suppose $Z := X \cap Y\neq \0$. Let $N$ be the set of all
vertices outside $Z$ with at least one neighbor in $Z$.
Let $\bf S$ be the (random) set of those
vertices of $N$ connected to $Z$ by at least one open edge.
We have
$$\Pr(A \, R_X) = \sum_S \Pr({\bf S} = S) \, \Pr(A \, R_X | S),$$
where the sum is over $S\sub N$ and we write
$\Pr( \cdot | S)$ for
$\Pr(\cdot | {\bf S} = S)$. Similar expressions hold for the
other terms in \eqref{prop1}.
Moreover, clearly,
$$
\Pr({\bf S} = S) \, \Pr({\bf S} = T)
=\Pr({\bf S} = S \cap T) \, \Pr({\bf S} = S \cup T)
~~~\forall S, T \sub N.
$$
So according to the Ahlswede-Daykin (``Four Functions")
Theorem (\cite{AD} or e.g. \cite{Boll}),
\eqref{prop1} will follow if we show that,
for all $S, T \sub N$,

\begin{equation}\label{suff}
\Pr(A \, R_X | S) \, \, \Pr(B \, R_Y | T) \leq
\Pr(A \, B \, R_{X \cap Y} | S \cap T) \,\,
\Pr(R_{X \cup Y} | S \cup T).
\end{equation}

Now it is easy to see that, for any set of vertices
$W \supseteq Z$,
and any event $D$
that does not depend on $E_W$,

\begin{equation}
\label{rest}
\mbox{$\Pr(D \, R_W | S) = \Pr'(D \, R_{(W \setminus Z) \cup S})$,}
\end{equation}
where $\Pr'$ refers to the induced model on the graph
$G'$ obtained from $G$ by removing $Z$.
(Strictly speaking, the $D$ on the r.h.s. of (\ref{rest})
is not the
same as that on the l.h.s., since it is a subset of
$\{0,1\}^{E \setminus E_Z}$ rather than $\{0,1\}^E$;
but since $D$ does not depend
on $E_W$ (and hence not on $E_Z$), the two events are essentially the same,
so we ignore the irrelevant distinction.)

Applying (\ref{rest}) to each of the four terms in (\ref{suff}),
we have

\begin{eqnarray*}
\Pr(A  R_X | S)  \Pr(B  R_Y | T) &=&
\mbox{$\Pr'(A  R_{(X \setminus Z) \cup S})
\Pr'(B  R_{(Y \setminus Z) \cup T})$}\\
&\leq& \mbox{
$\Pr'(ABR_{((X \setminus Z) \cup S) \cap ((Y \setminus Z) \cup T)})
\Pr'(R_{(X \setminus Z) \cup S \cup (Y \setminus Z) \cup T})$}\\
&\leq&
\mbox{$\Pr'(A  B  R_{((X \cap Y) \setminus Z) \cup (S \cap T)})
\Pr'(R_{((X \cup Y) \setminus Z) \cup (S \cup T)})$}\\
&=&\Pr(A  B  R_{X \cap Y} | S \cap T)
\Pr(R_{X \cup Y} | S \cup T),
%\label{suff2}
\end{eqnarray*}
where the first inequality follows from our inductive hypothesis
(applicable since $G'$ has fewer vertices than $G$),
and the second from

$$((X \setminus Z) \cup S) \cap ((Y \setminus Z) \cup T) \supseteq
((X \cap Y) \setminus Z) \cup (S \cap T)$$
and
$(X \setminus Z) \cup S \cup (Y \setminus Z) \cup T =
((X \cup Y) \setminus Z) \cup (S \cup T).$
\qed

\medskip
In particular we have the promised generalization of (\ref{vdBK1}):
\begin{thm}\label{1.2}
For s,A,B and $X$ as in Theorem~{\rm \ref{1.1}},
\begin{equation}
\label{thm1}
\Pr(A \, B \, | \, s \nrar X) ~ \geq ~ \Pr(A \, | \, s \nrar X) \,
\Pr(B \, | \, s \nrar X).
\end{equation}
\end{thm}
\begin{proof}
Take $Y = X $ in Theorem~\ref{1.1}.
\end{proof}
%It may be worth emphasizing that we don't know how to prove
%Theorem~\ref{1.2} without going through Theorem~\ref{1.1};
%indeed, generalizing in this way was the key idea in
%the original proof of (\ref{vdBK1}).
\noindent
{\bf Remarks}

1.  It is easy to see that Theorem \ref{1.2} is 
equivalent to the special case where $|X|=1$.  
(To reduce to this, simply identify
the vertices of $X$, retaining all edges connecting them to
$V\sm X$ (edges internal to $X$ may be deleted, but are anyway
irrelevant).)  We have used the present form both because 
it will be convenient for the proof of Theorem~\ref{1.5}, and
because it is natural from the point of view of the 
contact process application in Section~\ref{DPCP}.
Similarly, we could replace $s$ in all results of this section,
and $t$ in Theorems~\ref{1.4} and \ref{1.5}, by {\em sets} of vertices.
The same easy equivalence holds
for the directed graph results of Section~\ref{DPCP};
but in the case of the random cluster measures of Section~\ref{RCM}
the more general statements, while still true, do not seem to follow
in the same way from their specializations.

2.  The derivation of Theorem \ref{1.2}
may give the impression that it
is less general than Theorem \ref{1.1}, but in fact the two are equivalent.
To see this, first note that
(consideration of appropriate complementary events shows that)
Theorem \ref{1.2} also holds when $A$ and $B$ are both decreasing
rather than increasing, while the inequality (\ref{thm1})
reverses if one of $A, B$ is increasing and the other decreasing.
Thus for $A, B, X, Y$ as in Theorem \ref{1.1},
Theorem \ref{1.2} implies that conditioned
on $R_{X\cap Y}$, each of the pairs
$(A, R_{X\sm Y})$, $(B, R_{Y\sm X})$ is negatively correlated,
while each of $(A,B)$, $(R_{X\sm Y},R_{Y\sm X})$ is
positively correlated.  So, writing $\Pr'$ for our percolation
measure conditioned on $R_{X\cap Y}$, we have 
(compare (\ref{base}))

\begin{eqnarray*}
\mbox{$\Pr'(A \, R_{X \setminus Y} )
\Pr'(B \, R_{Y \setminus X} )$}
& \leq &
\mbox{$\Pr'(A)\Pr' (R_{X \setminus Y} )
\Pr'(B )\Pr'(R_{Y \setminus X} )$}\\
&\leq & \mbox{$\Pr'(A B )
\Pr'( R_{X \setminus Y}R_{Y \setminus X} ),$}
\end{eqnarray*}
which is equivalent to
\eqref{prop1}.

\bigskip
It will be helpful to have the ``functional extension"
of Theorem~\ref{1.2}:

\begin{thm}\label{TC}
For s, $X$ as in Theorem~{\rm \ref{1.2}}, and f,g
bounded, increasing, measurable functions of $C_s$,

$$
\E[f \, g \, | \, s \nrar X] \geq
\E[f \, | \, s \nrar X] \,
\E[g \, | \, s \nrar X].
$$

\mn
The inequality is reversed if one of f,g is increasing 
and the other decreasing.
\end{thm}

\begin{proof}  This is a standard (and easy) reduction.
We omit the argument for (a) and note that (b) is (a)
applied to the pair $(f,-g)$.
\end{proof}

\noindent

Recall that a collection of random variables
$\{\sigma_i: i \in \Ga\}$, with $\Ga$ a finite
or countably infinite index set and the $\sigma_i$'s
taking values in $\{0,1\}$ (or some other
ordered set),
is said to be {\em positively associated}
if for any two functions $f$, $g$ of the $\si_i$'s
that are either both increasing or both decreasing
(and, in case $\Ga$ is infinite, measurable),
one has $\E fg \geq \E f\E g$. The simplest non-trivial
example is when the $\sigma_i$'s are independent (Harris' inequality).

If we define a random subset $W $ of a set $T$
to be positively associated
if the collection $\{\eta(a) = {\bf 1}_{\{a \in W\}}: \, a \in T\}$
is positively associated, then
%If we define, for $a \in V$, $\eta(a):={\bf 1}_{\{s \rar a\}}$,
Theorem \ref{TC}
says that the open cluster of $s$ is
%the random variables
%$\eta(a), \, a \in V \setminus (X \cup \{s\})$, are
conditionally positively associated
given the event $\{s\nrar X\}$.
We will see further, similar examples later.

Positive association is often derived from the
FKG Inequality, which generalizes Harris' inequality and says that
positive association holds
for measures (on $\{0,1\}^n$, say) satisfying the
``positive lattice condition" (also called ``FKG lattice condition"), {\em viz.}
\begin{equation}
\label{plc}
\mu(\si)\mu(\tau)\leq \mu(\si\wedge \tau)\mu(\si\vee \tau)
\end{equation}
(where $(\si\wedge \tau)_x$ and $(\si\vee \tau)_x$ are the
minimum and maximum of $\si_x$ and $\tau_x$).
The positive lattice condition is much stronger than positive association.
It says that the conditional
probability that $\sigma_x = 1$, given the values
of $\sigma_y, y \neq x$, is increasing in those values.

%A point which may be of some interest here is that we
%are getting positive association in situations
%where (\ref{plc}) fails.
%Of course it is not hard to invent situations
%where this happens; but we are not aware of many natural
%examples of this type.
%{\bf [Rob:  I wasn't sure about this remark.  Feel
%free to modify or just delete.]}

Let us also
recall here that for measures $\nu$ and $\nu'$ on $\{0,1\}^n$ (or some other
partially ordered set), $\nu$ {\em stochastically dominates}
$\nu'$ ($\nu\succ \nu'$)
if $\nu(f) \geq \nu'(f)$ for every increasing function $f$
(where $\nu(f)$ is the expectation of $f$ w.r.t. $\nu$).

\bigskip
As suggested earlier, we do not see any good reason to expect
something like Theorem~\ref{1.2}.
(For instance, as noted in \cite{BeKa}, it is easy to see
that the analogous statement with $s\rar t$ in place of
$s\nrar t$ is false.)
Nonetheless, it implies the
following intuitively more natural statement,
%(and the aforementioned counterpart
%(\ref{new1}) of (\ref{vdBK1})),
which says, informally, that
conditioned on nonexistence of an open $(s,t)$-path,
the clusters
$C_s$ and $C_t$ are negatively correlated.

\begin{thm}\label{1.4}
Let $s$ and $t$ be (distinct) vertices, and $f$ and $g$
bounded measurable increasing functions of $C_s$ and $C_t$
respectively.  Then
$$
\E[f \, g \, | \, s \nrar t] \leq \E[f \, | \, s \nrar t] \,
\E[g \, | \, s \nrar t].
$$
\end{thm}

\noindent
Note that  \eqref{new1} is the special case where $f$ is the indicator
of the event
$\{s \rar a \}$, and $g$ that of the event $\{t \rar b \}$.

We have stated Theorem~\ref{1.4} above largely because, as mentioned
earlier, it was the original motivation for this work;
but the next statement, which contains Theorem~\ref{TC} 
as well as
Theorem~\ref{1.4}, seems to be
the correct level of generality here.

\begin{thm}\label{1.5}
Let $s$ and $t$ be (distinct) vertices, and $f$ and $g$
bounded, measurable functions of $(C_s,C_t)$, each increasing
in $C_s$ and decreasing in $C_t$.  Then
\begin{equation}
\label{X}
\E[f \, g \, | \, s \nrar t] \geq \E[f \, | \, s \nrar t] \,
\E[g \, | \, s \nrar t].
\end{equation}
\end{thm}
\noindent
In other words, on $\{s\nrar t\}$ we have positive 
association of all the r.v.'s ${\bf 1}_{\{e\in C_s\}}$ and 
${\bf 1}_{\{e\not\in C_t\}}$.
(Note: here and often in what follows, we use
``on $Q$" to mean ``conditioned on (the event) $Q$.")

\mn
{\em Proof of Theorem~\ref{1.5}.}
As in the case of Theorem 1.1, it is enough to prove
this for finite $G$.

%It is convenient to work with the (equivalent) functional
%version; that is, to show that for $f,g$ increasing functions
%of $C_s$ and $C_t$ respectively,
%\begin{equation}\label{X}
%\E[f \, g \, | \, s \nrar t] ~\leq ~
%\E[f \, | \, s \nrar t] \,\, \E[g \, | \, s \nrar t]
%\end{equation}
%(which of course contains (\ref{thm2}) as the case
%$f={\bf 1}_A,~g={\bf 1}_B$).

We have
\begin{equation}\label{expand}
\E[f \, g\, | \, s \nrar t]
= \sum_W \Pr(C_s = W \, | \, s \nrar t) \,\, \E[ f \, g\, | \, C_s = W],
\end{equation}
where we may restrict to $W$ containing no $(s,t)$-path.
Write $\bar{W}$ for the union of $W$ and its ``boundary";
that is, $\bar{W}$ consists of all edges
having
at least one vertex in common with some edge of $W$.

When we condition on $\{C_s = W\}$, $f$ and $g$ become decreasing
functions of $C_t$, and the (conditional) distribution of $C_t$
is the same as that for the restriction of our percolation model
to the graph obtained from $G$ by deleting all edges in $\bar{W}$.
Thus (on $\{C_s = W\}$) $f,g$ are decreasing functions of the
independent r.v.'s $(\om_e:e\in E\sm \bar{W})$,
and by Harris' inequality we have
$\E[fg|C_s=W]\geq \E[f|C_s=W]\E[g|C_s=W]$.

On the other hand, the conditional distribution of
$C_t$ given $\{C_s = W\}$ is stochastically decreasing in
$W$ (to couple these distributions, choose {\em all}
$\om_e$'s independently according to their $p_e$'s and
then for conditioning on $\{C_s = W\}$ simply ignore
those $\om_e$'s with $e\in \bar{W}$);
so in particular
$\E[f|C_s=W]$ and $\E[g|C_s=W]$ are increasing functions of $W$,
and it then follows from
Theorem~\ref{TC} that the right hand side of
(\ref{expand}) is not less than

$$
\left(\sum_W \Pr(C_s = W \, | \, s \nrar t) \,
\E[f|C_s=W] \right)~
\left(\sum_W [\Pr(C_s = W \, | \, s \nrar t) \,
\E[g|C_s=W]\right)~~~~~~~~~
$$
$$
~~~~~~~~~~~~~~~~~~~~~
= \E[f \, | \, s \nrar t] \,\, \E[g \, | \, s \nrar t];
$$
so we have \eqref{X}.
\qed

As just shown, Theorem \ref{1.4} follows easily 
from Theorem \ref{TC}.
While one might expect a similar proof (or {\em some} proof)
of the reverse implication,
we do not see this.
In Section 2 we will (as mentioned earlier) take a completely different
approach which, even for the more general class of random-cluster
measures, gives Theorems \ref{TC} and \ref{1.4} ``simultaneously".

\section{Random-cluster measures}
\label{RCM}
\subsection{Definitions and a Markov chain proof}\label{DMCP}
A well-known generalization of the bond percolation model is the
{\em random-cluster model} (RCM) introduced by Fortuin and Kasteleyn circa 1970.
(See e.g. \cite{Grimmett}, Section 13.6, \cite{Grimmett3}
for additional background and references.)

Let $G = (V,E)$ be a finite graph.
In addition to the parameters
$p_e, \, e \in E$ of the
ordinary bond percolation model, the random-cluster model is equipped
with a positive parameter $q$. To avoid trivialities we assume that
$0 < p_e < 1$ for all $e \in E$.
The {\em random-cluster measure} (r.c.m.)
with the above parameters on $\Om = \{0,1\}^E$ is then given by

\begin{equation}
\label{rcdef}
\varphi_q(\om) \, (=\varphi_{G,q}(\om)) \, \propto
\,  q^{k(\om)}  \prod_{e \in E \, : \, \om_e = 1} p_e
\prod_{e \in E\, : \, \om_e = 0} (1-p_e), \,\,\,\, \,\,
\om \in \Om,
\end{equation}

\noindent
where $k(\om)$ is the number
of connected components in $\om$, and,
as usual, $f(\om) \propto g(\om)$ means
$f(\om) = C g(\om)$ for some (positive) constant $C$.
(For the present discussion we regard the $p_e$'s as
given once and for all, and omit them from our notation.)

Thus $q = 1$ gives the ordinary bond
percolation model.
We have, in spite of serious attempts, 
not been able to adapt the approach
of Section 1 to $q>1$.
(We do not consider $q<1$, for which the correlation 
properties of the model are quite different).
Here we take a different, 
``dynamical" approach, based on the introduction of a 
Markov chain whose states are pairs of clusters
(this is not the only possibility;
see the remark following the proof of Theorem~\ref{2.5})
which converges to a measure (on pairs of clusters) 
corresponding to (\ref{rcdef}), 
and for intermediate stages 
of which the correlation properties we are after can be derived
%via a simple (presumably not new) lemma on positive association
%(Lemma~\ref{LPA})
from known properties of the RCM.

For the following extension of Theorem~\ref{1.5} to the RCM
we replace the vertices $s$ and $t$ by sets $S$ and $T$, 
recalling that
the remark following Theorem~\ref{1.2} regarding
the easy reduction from sets to singletons
is not valid here.
Extending our earlier notation, 
we use $C_S$ for the set of
edges belonging to open paths
starting at vertices of $S$.

\begin{thm}\label{2.5}
Consider a distribution {\rm (\ref{rcdef})} with $q\geq 1$.
Let $S$ and $T$ be disjoint sets of vertices, and $f$ and $g$
bounded, measurable functions of $(C_S,C_T)$, each increasing
in $C_S$ and decreasing in $C_T$.  Then on $\{S\nrar T\}$,
\begin{equation}
\label{X'}
\E f g \geq \E f \E g.
\end{equation}
\end{thm}

Following the Markov chain proof of this, we also
give, in Section~\ref{Fuzzy}, a different
argument, which unfortunately seems only to work when $q \geq 2$
and $|S|=|T|=1$.
So, somewhat strangely, we have separate (and distinct) proofs
for the cases $q =1$ and
$q \geq 2$, but for the intermediate
case $1 < q < 2$, no alternative to the Markov chain approach.

\begin{proof}
We first give some additional notation, and state
some (well-known) lemmas we will need.
If $F$ is a subset of $E$ (the set of edges of our graph $G$), and $\om \in
\Om = \{0,1\}^E$, we write $\om_F$ for the restriction of $\omega$ to $F$
($\om_F=(\om_e: e \in F)$), and $V(F)$ for the set of vertices 
incident with edges of $F$.
We continue to use the notation $\bar{W}$ 
introduced following \eqref{expand}.

\begin{lem}\label{lem3.1}
For $q \geq 1$, the random-cluster measure \eqref{rcdef}
satisfies the positive lattice condition \eqref{plc}.
\end{lem}

When $\phi_{G,q}$ is conditioned on the values of
some of the variables $\om_e$, 
the remaining variables are distributed as they 
would be under the (natural) 
r.c.m. on the graph obtained from $G$
by deleting $e$'s with $\om_e=0$ and contracting those
with $\om_e=1$.  For our purposes the relevant 
cases of this are given by

\begin{lem}\label{3.2alt}
Let $A \subset V$ and $F \subset E$.
The restriction of 
$\varphi_{G,q}$ to $\{0,1\}^{E \setminus \bar{F}}$
under conditioning on either of the events $\{C_A = F\}$,
$\{\om_{\bar{F}}\equiv 0 \}$
(i.e. $\{\om_e =0\forall e\in\bar{F}\}$)
is the r.c.m. with parameter $q$ on $G -\bar{F}$ 
(the graph obtained 
from $G$ by deleting all edges in $\bar{F}$); more formally,
$$
\phi_{G,q}(\om_{E \setminus \bar{F}} 
= \cdot \, | \, C_A = F) =
\phi_{G,q}(\om_{E \setminus \bar{F}} 
= \cdot \, | \, \om_{\bar{F}}\equiv 0) =
\varphi_{G -\bar{F}, q}(\cdot).
$$
\end{lem}
\mn
If $A,F$ are as in Lemma~\ref{3.2alt}, and $B\sub V\sm V(F)$,
then $\{C_A=F\}\sub \{A\nrar B\}$; so Lemma~\ref{3.2alt}
implies

\begin{lem}\label{3.3alt}
If $A$, $F$ are as in Lemma \ref{3.2alt},
$B\sub V\sm V(F)$, and $\phi$ is (temporarily)
$\phi_{G,q}$ conditioned on $\{A\nrar B\}$,
then
$$
\phi(\om_{E \setminus \bar{F}} 
= \cdot \, | \, C_A = F) =
\phi_{G,q}(\om_{E \setminus \bar{F}} 
= \cdot \, | \, \om_{\bar{F}}\equiv 0) =
\varphi_{G -\bar{F}, q}(\cdot).
$$
\end{lem}

\medskip\noindent
We now turn to the proof of Theorem~\ref{2.5}.
We consider a Markov chain with state space
$\hOm$ consisting of pairs $(C_S,C_T)$ satisfying $Q:=\{S\nrar T\}$.
(So the states are pairs $(C,C')$ such that
$C,C'\sub E$; $C$ (resp. $C'$)
is a union of paths beginning at vertices of $S$
(resp. vertices of $T$); and $V(C)\cap V(C')=\0$.)

We write $\phi$ for the 
measure $\phi_{G,q}$ conditioned on $Q$ and
$\hap$ for the measure which $\phi$ induces on $\hOm$.
%We will, however, abuse this and related notation slightly:
%for an event $A$
%depending on $\om\in\{0,1\}^E$ we will sometimes use
%``$\hap$ conditioned on $A$" as shorthand for
%``the measure which $\phi$ conditioned on $A$ 
%induces on $\hOm$."

Initially our chain is in some fixed state $(C_S^0,C_T^0)\in\hOm$.
Given $(C_S^{i-1},C_T^{i-1})$, the state of the chain at time
$i-1$, we choose $(C_S^i,C_T^i)$
in two steps, first choosing $C_T^i$ according to $\phi$ conditioned
on $\{C_S=C_S^{i-1}\}$---that is,
$$
\Pr(C_T^i= \cdot) =\phi(C_T = \cdot|C_S=C_S^{i-1})
$$
---and then, similarly, $C_S^i$ according to
$$
\Pr(C_S^i= \cdot) =\phi(C_S = \cdot|C_T=C_T^i).
$$
It is clear that $\hap$ is stationary for this chain, and that the
chain is irreducible and aperiodic;
so to prove Theorem~\ref{2.5} it's enough to show

\begin{claim}\label{cl}
For $f,g$ as in the statement of Theorem~\ref{2.5} and any n,
{\rm (\ref{X'})} holds for expectation taken with respect to the law of
$(C_S^n,C_T^n)$.
\end{claim}

\mn
Let $X_e^i$ and $Y_e^i$ be the indicators of the events 
$\{e\not\in C_T^i\}$
and $\{e\in C_S^i\}$ ($e\in E$, $i=0,1,\ldots$).
These are, of course, not independent, 
but we will show, using the following
presumably well-known observation, that they are 
positively associated.

\begin{lem}\label{LPA}
Suppose $W_1\dots W_a$ and $Z_1\dots Z_b$ are (say) 
$\{0,1\}$-valued r.v.'s
with joint distribution $\psi$ satisfying

\mn
{\rm (i)}  $W_1\dots W_a$ are positively associated;

\mn
{\rm (ii)}  $Z_1\dots Z_b$ are conditionally positively 
associated given
$W_1\dots W_a$; and

\mn
{\rm (iii)}  for $W,W'\in\{0,1\}^a$ with $W'\geq W$,
$\psi(\cdot|W')\succ \psi(\cdot | W)$,
where $\psi(\cdot |W)$ is the conditional distribution of 
$(Z_1\dots Z_b)$
given $(W_1\dots W_a) = W$.

\mn
Then
$W_1\dots W_a,Z_1\dots Z_b$ are
positively associated.
\end{lem}
\begin{proof}
Suppose $f,g$ are increasing functions of $W_1\dots Z_b$,
and for $W\in\{0,1\}^a$, set
$F(W)=\E [f|W]$ ($:= \E[f|(W_1\dots W_a)=W]$)
and $G(W)=\E [g|W]$.  Then
\begin{eqnarray*}
\E fg &=& \E\{\E[fg|W]\}\\
&\geq & \E\{\E[f|W]\E[g|W]\}\\
&\geq & \E F\E G \\
&=& \E f\E g,
\end{eqnarray*}
where the first inequality follows from (ii) and the second
from (iii) and (i).
\end{proof}

\begin{lem}\label{PA}
The collection
\begin{equation}\label{coll}
\cup_{i\geq 1}\cup_{e\in E}\{X_e^i,Y_e^i\}
\end{equation}
is positively associated.
\end{lem}
Note this is enough for Claim~\ref{cl} since (trivially)

\begin{rem}\label{inc}
For each $n$, $C_S^n$ is increasing in the variables
$X_e^i$, $Y_e^i$, and $C_T^n$ is decreasing in these variables.
\end{rem}

\noindent
{\em Proof of Lemma~\ref{PA}.}
Of course it's enough to show positive association for finite 
subsets of the collection
(\ref{coll}).
We will show by induction on $n$ that for
each $n$, each of the collections
\begin{equation}\label{coll1}
\{X_e^i:e\in E, i\leq n\}\cup \{Y_e^i:e\in E, i<n\}
\end{equation}
and
\begin{equation}\label{coll2}
\{X_e^i:e\in E, i\leq n\}\cup \{Y_e^i:e\in E, i\leq n\}
\end{equation}
is positively associated.
(The base cases---those with $n=0$---are, of course, trivial.)
Actually we just give the argument for (\ref{coll1}), that for
(\ref{coll2}) being essentially the same.

We want to apply Lemma~\ref{LPA} with
$(W_1\dots W_a)=\cup((X_e^i,Y_e^i):e\in E, i<n)$
and $(Z_1\dots Z_b)=(X_e^n:e\in E)$,
so need to verify conditions (i)-(iii) of the lemma.
Of course (i) is just our inductive hypothesis, so
our concern is really with (ii) and (iii).
%we appeal to what we know about 
%random cluster measures.

Consider a possible value $W$ of $(W_1\dots W_a)$,
with $F$ the corresponding value of $C_S^{n-1}$.
Under conditioning on $\{(W_1\dots W_a)=W\}$, we 
have $X_e^n$ fixed for $e\in \bar{F}$ 
(namely $X_e^n\equiv 1 ~\forall e\in \bar{F}$), 
while, by Lemma~\ref{3.3alt}, 
the remaining $X_e^n$'s are distributed as the 
variables ${\bf 1}_{\{e\not\in C_t(\om)\}}$, 
where $(\om_e:e\in E\sm \bar{F})$ 
is chosen according to $\phi_{G-\bar{F},q}$.
Positive association of these variables is
given by Lemma~\ref{lem3.1}, so we have (ii).

Now let $W'$ be a second possible value of $(W_1\dots W_a)$,
with $W'\geq W$ and $F'$ the corresponding value of $C_S^{n-1}$.
According to Remark~\ref{inc} 
we have $F\sub F'$.
So (iii) amounts to saying that for $F\sub F'\sub E$ and 
$h$ any increasing function of $C_T$,
\begin{equation}\label{FF'}
\phi(h|C_S=F)\geq \phi(h|C_S=F')
\end{equation}
(note $h$ is a {\em decreasing} function of the $X_e^n$'s).
But using Lemmas~\ref{3.2alt} and \ref{3.3alt}, 
we may rewrite the left and right hand sides of (\ref{FF'})
as
$$
\phi_{G-\bar{F},q}(h)
$$
and 
$$
\phi_{G-\bar{F'},q}(h)=
\phi_{G-\bar{F},q}(h|\mbox{$\om_{\bar{F'}\sm\bar{F}}\equiv 0$});
$$
and then (\ref{FF'}) follows from Lemma~\ref{lem3.1}
(which gives positive association for the measure 
$\phi_{G-\bar{F},q}$).\end{proof}

\mn
{\bf Remark.}
We briefly indicate an alternative proof
of Theorem~\ref{2.5}, again using a Markov 
chain and based on a similar idea.
This, our original proof, 
is perhaps more natural than that given 
above, but does not seem
as easily adapted to prove the directed version
of Theorem~\ref{1.5} (Theorem~\ref{3.5}).

We again use $\phi$ for $\phi_{G,q}$ conditioned on 
$\{S\nrar T\}$.
Our chain in this case is $\om^0, \om^1, \om^2, \ldots$
drawn from the state space 
$\hat{\Om} := \{\om \in \Om \,:\, S \nrar T\}$.
Initially the chain is in some fixed state $\om^0$.
Given $\om^{i-1}$, the state of the chain at time
$i-1$, we choose $\om^i$
in two steps, first choosing an intermediate configuration $\tau^i$
according to $\phi$ conditioned
on $\{C_S=C_S(\om^{i-1})\}$---that is, for $\ze\in\hat{\Om}$
with $C_S(\ze)=C_S(\om^{i-1})$,
$$
\Pr(\tau^i= \ze) =\phi(\om = \ze|C_S(\om) =C_S(\om^{i-1}))
$$
---and then, similarly, $\om^i$ according to
$$
\Pr(\om^i= \ze) =\phi(\om = \ze|C_T(\om)=C_T(\tau^i)).
$$
It is clear that $\phi$ is stationary for this chain, and that the
chain is irreducible and aperiodic;
so to prove Theorem~\ref{2.5} it's enough to show

\bn
{\bf Claim.}
{\em For $f,g$ as in the statement of 
Theorem~\ref{2.5} and any n,
{\rm (\ref{X'})} holds for expectation taken with 
respect to the law of $(\om^n)$.}

\bn
To prove this we introduce independent r.v.'s
$X_e^i$, $Y_e^i$ ($e\in E$, $i=1,\ldots$),
each uniform on $[0,1]$, and some fixed ordering,
``$\prec$," of $E$.
Then to decide the value of $\tau^i_e$ we 
compute the conditional probability, say $\alpha$,
that $\tau^i_e =1$ given the values of the
$\om_e^{i-1}$'s (or just the value of $C_S(\om^{i-1})$)
and those $\tau_{e'}^i$'s with
$e'\prec e$, and set $\tau_e^i=1$ iff $X_e^i<\alpha$.
For $\om^i$ we proceed analogously, with the
requirement for $\om_e^i=1$ now being $Y_e^i>1-\alpha$.

It is then not hard to show, again using 
Lemmas~\ref{lem3.1}-\ref{3.3alt}, that (for each $n$)
$\om^n$ is increasing in the variables $X_e^i$, $Y_e^i$,
so that the Claim follows from Harris' inequality.

\subsection{A separate proof for $q \geq 2$}\label{Fuzzy}
As mentioned earlier, it turns out, somewhat curiously,
that for $q\geq 2$ and $S$ and $T$ consisting of 
single vertices $s$ and $t$,
we can prove Theorem~\ref{2.5} in a different
way by
exploiting a connection between the random cluster
model and the ``fuzzy Potts model."
(The corresponding
connection involving the ordinary Potts model again
goes back to Fortuin and Kasteleyn.) 
Before doing so, we need to
review some classical and more
recent facts concerning this connection.

\smallskip\noindent
Let $q=\al+\be$ with $\al,\be>0$. Using the random-cluster measure $\varphi$
we generate a random spin configuration $\sigma \in \{0,1\}^V$ 
as follows.

\mn
(i)  Choose $\omega \in \{0,1\}^E$ according to
$\varphi_q$.

\mn
(ii)  For each component $C$ of $\om$, let $\si$ take
the value 1 (resp. 0) on all vertices of $C$ with
probability $\al/q$ (resp. $\be/q$), independently of
the values of $\si$ on other components.
Let $\mu_{\al,\be}$ denote the distribution of $\si$.

\mn
(In \cite{Haggstrom} this is called the {\em fractional
fuzzy Potts model}.)
This procedure produces a coupling measure 
$\cal P$ of $\om$ and $\si$,
or, rather, of $\varphi_q$ and $\mu_{\al,\be}$.
So we may also think of first choosing $\si$ and then drawing
from the conditional distribution 
$\cal P(\cdot \, | \, \si)$ to obtain
a typical (with distribution $\varphi_q$) 
edge configuration $\om$.
It is known
(and easy to check) that this ``reversed" 
procedure can be described as follows.

\mn
(iii)  Choose $\si$ according to $\mu_{\al,\be}$.

\mn
(iv)  For $i=1,0$, let $G(i)=G[\si^{-1}(i)]$
(the ({\em induced}) subgraph consisting of vertices in
$\si^{-1}(i)$ and edges of $G$ contained in 
this set).
Set $\om_e=0$ whenever $\si$ assigns different
values to the ends of $e\in E$, and choose the 
restrictions of $\om$ to $E(G(1))$ and $E(G(0))$
(independently) according to 
$\varphi_{G(1),\al}$ and $\varphi_{G(0),\be}$.\\

%\mn
%(iv) The edge set $E$ is clearly a partition of the sets $E(1)$
%(the set of those edges
%of which both endpoints have $\si$-value $1$), $E(0)$ (similar,
%but now for value $0$) and $E(01)$
%(the set of edges whose two endpoints have different
% $\si$-values). Let $G(1)$ be the graph with edges $E(1)$ 
%and vertices all
%endpoints of edges in $E(1)$. Similarly define $G(0)$.
%Choose a random configuration $\om^{(1)} \in \{0,1\}^{E(1)}$
%according to $\varphi_{G(1),\al}$ and, independently, a random
%configuration $\om^{(0)} \in \{0,1\}^{E(0)}$
%according to $\varphi_{G(0),\be}$.
%Finally, the configuration $\om \in \{0,1\}^E$ is 
%now obtained by taking
%$\om_e$ equal to $0$ if $e \in E(01)$, to $\om^{(1)}_e$ if
%$e \in E(1)$ and to $\om^{(0)}_e$ if
%$e \in E(0)$. \\
Furthermore, if in (iii) we choose $\si$ according to the 
{\em conditional}
distribution
$$\hat{\mu}_{\al,\be}(\cdot) = 
\mu_{\al,\be}(\cdot \, | \,\si(s)=1,\si(t)=0), $$
then $\om$ (in (iv)) has the
distribution we want, namely
$\varphi_{G,q} (\cdot \, | \, s \nrar t)$;
so for Theorem~\ref{2.5} we may take 
$\om$ to be chosen in this way.

\medskip
The salient points for our purposes are
then as follows.
%(with positive association as
%defined below Theorem \ref{TC}). 
(For clarity
we now add subscripts to the expectation symbol
$\E$ to indicate measures with respect 
to which expectation is
taken.)

\mn
(a)  (Lemma \ref{lem3.1})
For any graph $H$ and $c\geq 1$, $\varphi_{H,c}$
has positive association.  

\mn
(b)  It is shown in \cite{Haggstrom} that for any
$\al,\be\geq 1$, $\mu=\mu_{\al,\be}$ satisfies
(\ref{plc}), whence, according to the FKG Inequality,
$\mu_{\al,\be}$ and
$\hat{\mu}_{\al,\be}$ are positively associated.
%the $\si_x$'s are positively associated and remain so under
%(\ref{Q}).

\mn
(c)  If $\al,\be\geq 1$ and
$f$ is a function of $(C_s,C_t)$ which is 
increasing in $C_s$ and decreasing in 
$C_t$, then
$\E_{\cal P}[f\, |\,\si]$
is increasing $\si$.
% and the same holds if $\cal P$ is replaced by $\hat{\cal P}$.
(This follows from ((iv) and) (a).)
%, the assertion for
%$f$ (resp. $g$) actually requiring only $\al\geq 1$
%(resp. $\be\geq 1$).)

\bn
{\em Alternate proof of Theorem~\ref{2.5} 
for $q\geq 2$, $S=\{s\}$ and $T=\{t\}$.}
%For Theorems 1.3 and 1.4 with $q\geq 2$
%we may apply this machinery with any $\al,\be\geq 1$;
%for instance, for Theorem 1.3:  for
%\begin{thm}\label{rcthm}
%Theorem~\ref{TC} and \ref{1.4} hold for the random-cluster distribution
%(\ref{rcdef}) with $q \geq 2$.
%\end{thm}

Let $f,g$ be as in the statement of the theorem.
Fix some $\alpha, \beta \geq 1$ with 
$\alpha + \beta = q$.
%(The most natural choice is $\alpha = \beta = q/2$. When 
%$q = 2$, where
%the corresponding Potts model is the Ising model, 
%this special choice
%makes the above machinery somewhat more elegant, but for
%general $q \geq 2$ there is no ``ideal" choice).
For simplicity we write $\mu$ for $\mu_{\al,\be}$ and
$\varphi$ for $\varphi_q$. The 
connections described above give

\begin{eqnarray*}
\E_{\varphi}[f g \, | \, s \nrar t] &=&
\sum_{\si}\hat{\mu}(\si)\, \E_{\cal P}[f(\om)g(\om) |  \si]\\
&\geq & 
\sum_{\si}\hat{\mu}(\si) \, \E_{\cal P}[f(\om) | \si] \,\,
\E_{\cal P}[g(\om) |  \si]\\
&\geq & 
\sum_{\si}\hat{\mu}(\si) \, \E_{\cal P}[f(\om) |  \si] \,
\sum_{\si}\hat{\mu}(\si) \, \E_{\cal P}[g(\om)|\si]\\
&=&
\E_{\varphi} [f \, | \, s \nrar t] 
  ~\E_{\varphi} [g  \, | \, s \nrar t],
\end{eqnarray*}
where the first inequality follows from (a) (and (iv))
and the second from (b) and (c).
\qed

\section{Directed percolation and contact processes}
\label{DPCP}
In this section we consider another generalization of ordinary
percolation: as in Section 1 we have a product distribution on
$\{0,1\}^E$, but now some (or all, or none) of the edges of our
graph are oriented. There are (at least) two natural ways to try
to extend the results of Section 1 to this setting, 
corresponding
to two possible extensions of the conditioning event
$\{s\nrar t\}$.
As we will see, both extensions are reasonable for
Theorems~\ref{1.1}-\ref{TC}, but only one of them makes 
sense for
Theorem~\ref{1.4} (and Theorem~\ref{1.5}).
The first set of extensions yield in particular improvements of
some of the results of Belitsky, Ferrari, Konno and Liggett 
\cite{BFKL}
regarding the contact process (defined below).
We will first indicate these extensions (proofs of which are essentially
identical to the proofs of the corresponding statements in Section~\ref{Intro}) and
discuss their relevance to the contact process, before 
turning to the second set of extensions.

\medskip
We will need the following additional notation.
Unoriented edges will be denoted by $\{v,w\}$
and oriented edges by $(v,w)$ (where the orientation is
from $v$ to $w$).
When we speak of a {\em path}, we will now mean one which respects the orientations
of its oriented edges.  We write
$\{s \ra t \}$ for the event that there is an open path from $s$ to $t$
and $\{s \nra t \}$ for the complement of this event.
The {\em open cluster}, $C_s$, of $s$ is again
the set of all edges contained in
open paths starting at $s$.
As in Section 1, we fix a vertex $s$ and set
$R_X =\{s \nra x \, \forall x \in X\}$
for each $X\sub V\sm\{s\}$.
Of course all these definitions collapse to those of
Section~\ref{Intro} in case there
are no oriented edges; so the next result contains 
Theorem~\ref{1.1}.

\begin{thm}\label{di1}
With the preceding modified definitions,
Theorem~\ref{1.1} holds for directed percolation.
\end{thm}

\mn
{\em Proof.}
The proof is essentially the same as that of Theorem 1.1,
the only difference being that we should now take $N$ to be the set of
those $i\not\in Z$ for which there is at least one edge $(i,j)$ or $\{i,j\}$
with $j\in Z$, and modify the definition of $S$ similarly.\qed

\medskip
Our first extension of Theorem \ref{TC} to directed 
percolation is

\begin{thm}\label{3.3}
Let $s \in V$, $X \subseteq V \setminus \{s\}$ , and $f,g$ increasing
functions of $C_s$.
Then on $\{s\nra X\}$,
$$
\E fg \geq \E f\E g.
$$
\end{thm}
\noindent
This can be derived from Theorem~\ref{di1} in the same way 
as Theorem~\ref{TC} was derived from Theorem~\ref{1.1}.
It can also be proved using Markov chains 
(following either the proof of Theorem~\ref{2.5} or 
the alternate sketched afterwards), where we should
now replace $C_t$ (in its various incarnations) 
by the set of edges in paths {\em ending in t}.

\medskip\noindent
{\bf Remarks, and consequences for Contact Processes}

\smallskip\noindent
(i)
Analogously to what we said in Section 1,
Theorem \ref{3.3} can be stated in terms of
(conditional)
positive association; namely for any $X\sub V$,
the random variables
$\eta(y):= {\bf 1}_{\{s \ra y\}}$,
conditioned on the event $\{\eta \equiv 0 $ on $X\}$, are postively associated.

\smallskip\noindent
(ii) Taking $A = B = \Om$ in Theorem \ref{di1} gives

\begin{equation}
\label{conv}
\Pr(R_X) \, \Pr(R_Y) \leq \Pr(R_{X \cup Y}) \, \Pr(R_{X \cap Y}).
\end{equation}
%(Alternatively, though more clumsily, this follows from
%the aforementioned positive association:
%if $X',Y' \subset V \setminus \{s\}$ and we take $A=\ov{R}_{X'\sm Y'}$ and
%$B=\ov{R}_{Y'\sm X'}$  and $X = Y = X' \cap Y'$,
%then an instance of Theorem 2.1 is
%$$
%\Pr(AR_{X'\cap Y'})\Pr(BR_{X'\cap Y'})
%\leq \Pr(ABR_{X'\cap Y'})\Pr(R_{X'\cap Y'}).
%$$
%This last inequality clearly
%remains valid after replacing $A$ and $B$ by their complements,
%which gives another form of \eqref{conv}.)

\smallskip\noindent
(iii)
Belitsky, Ferrari, Konno and Liggett
(\cite{BFKL}, Theorem 1.5)
proved a special case of
\eqref{conv}
involving a particular graph on the vertex set $\Z^2 $.
Their argument actually applies whenever $V$ admits a
partition
$(V_0=\{s\})\cup V_1 \cup \cdots$ such that each
edge is directed from
$V_{i-1} $ to $V_i$ for some $i$ and $X\cup Y$ is
contained in some $V_i$, but does seem to depend
essentially on these properties.

\smallskip\noindent
(iv)
Much of \cite{BFKL} deals with the
% so-called
{\em contact process}
on a countable set $S$.
%a simply
%formulated---but mathematically
%very rich---spatial/probabilistic abstraction of the
%evolution of an infection in a population.
See \cite{Liggett} and \cite{Liggettnew} for background on this model;
very briefly:
Each site (individual) in $S$
can be in either of the states $1$
(ill and contagious) or $0$ (healthy, noncontagious).
Time is continuous, with $\eta_t(x)$ denoting
the state of site $x$ at time $t$.
An infected site $x$ becomes
healthy at rate
$\de_x$, and
a healthy site becomes ill at rate
$\sum_y \lambda(x, y) \eta(y)$.
Here $\de_x, x \in S$ and $\lambda(x, y), \, x, y \in S$, are the
parameters of the model. They are assumed to be non-negative and,
if $S$ is infinite, to satisfy the following
conditions (see \cite{BFKL}):
$\sup_{x \in S} \delta(x) < \infty$ and $\sup_{x \in S} \sum_{y \in S}
[\lambda(x,y) + \lambda(y,x)] < \infty$. 

A nice aspect of the model is that it can be viewed
in terms of percolation,
via a graphical representation (see e.g. \cite{Liggettnew}, pages 32-34):
being ill at some given time
corresponds to the existence of an appropriate path
in space-time. In fact, as is well-known,
the process can, by time-discretization, be
approximated by a directed percolation model.
(See the subsection on correlation inequalities, in particular
page 11, of \cite{Liggettnew} for the general idea of how
correlation inequalities for collections of independent Bernoulli
random variables can be
extended to continuous-time interacting particle systems,
and page 65 of \cite{Liggettnew} for a
concrete example for the contact process).
Combining this with the present results,
one obtains, in a straightforward way,
contact process analogues of the conditional
association property stated in (ii) above.
In particular this gives the following theorem.
\begin{thm}\label{contact}
Suppose $(\eta_t: t \geq 0)$ is a contact process as above,
with
deterministic initial configuration $\eta_0$.
Then for each $W \subset S$ and $t \geq 0$,
the collection $(\eta_t(x): x \in S \setminus W)$ is
conditionally positively associated given
$\{\eta_t \equiv 0$ on $W\}$.
\end{thm}
(An example of Liggett \cite{Liggett94} shows that if we instead
condition on $\{\eta_t \equiv 1$ on $W\}$, the above positive
association need not hold.)

Suppose now that at time $0$ {\em each} site is ill.
Let $\nu_t$ be the law of $\eta_t$
($= (\eta_t(x): x \in S)$).
It is well-known
(and follows easily from standard monotonicity arguments)
that as $ t \ra \infty$,
$\nu_t$ tends to a limit,
called the {\em upper invariant measure}
of the process and denoted $\nu$.
Clearly the preceding conditional association property for finite times extends to $\nu$.
So, if $W \subset S$, and $A$ and $B$ are events that
are determined by, and that are both increasing (or both decreasing) in
the $\eta(x)$, $x \in S \setminus W$, then
\begin{equation}
\label{passoc}
\nu(A \, B \,| \,\eta\equiv 0\,\mbox{on}\, W)
\geq \nu(A \, | \,\eta\equiv 0\,\mbox{on}\, W) \,\,
\nu(B \, | \,\eta\equiv 0\,\mbox{on}\, W).
\end{equation}
This is
a considerable strengthening of
a conjecture of Konno (\cite{Konno}, Conjecture 3.4.13),
which was proved in---and seems to have been the main
motivation for---\cite{BFKL}
(see inequality (1.3) in \cite{BFKL}), namely:
for any $K,L\sub S$,
%eqref{passoc}:
%$$ \nu(\eta_i = 0, i \in K \cup L \, | \, \eta_i = 0, \, i \in K \cap L)
%\geq \nu(\eta_i = 0, i \in K \, | \, \eta_i = 0, \, i \in K \cap L) \,\,
%\nu(\eta_i = 0, i \in L \, | \, \eta_i = 0, \, i \in K \cap L)$$,
%or, equivalently, that
\begin{equation}
\label{bellcase}
\nu(K \cap L) \,\, \nu(K \cup L) \geq \nu(K) \, \nu(L),
\end{equation}
where, for $M \subseteq S$,
$\nu(M):=\nu\{\eta\,:\,\eta\equiv 0\,\mbox{on}\,  M\}$.
(Of course \eqref{bellcase} is the special case of
\eqref{passoc}
in which $W=K\cap L$,
$A=\{\,\eta\equiv 0\,\mbox{on}\, K \setminus L\}$ and
$B=\{\,\eta\equiv 0\,\mbox{on}\, L \setminus K\}$.)

%Note that we could also derive \eqref{bellcase} from Theorem 2.2 in a
%different way,
%namely via \eqref{conv}.
%We chose the present route via \eqref{passoc},  which we think is a
%natural and
%interesting-in-itself strengthening of the above mentioned result in
%\cite{BFKL}.

\bn
{\bf A directed version of Theorem~\ref{1.5}}

For a sensible generalization of Theorems~\ref{1.4} and \ref{1.5} to the present
setting we need a different substitute for $\{s\nrar t\}$.
It is easy to see that neither $\{s\nra t\}$ nor $\{s\nra t\nra s\}$
will do here (e.g.
consider the graph on $\{s,t,v,a\}$ with (oriented) edges $(s,v),(t,v),(v,a)$, and
events $A=\{s\ra a\}$, $B=\{t\ra a\}$);
but there is another natural choice which does work, at least when we assume
there are no undirected edges.
Recall $V(F)$ is the set of vertices incident with edges of $F$.

\begin{thm}\label{3.5}
Assume G is a digraph in the usual sense (that is, all its
edges are directed).
Let $s$ and $t$ be (distinct) vertices, and $f$ and $g$
bounded, measurable functions of $(C_s,C_t)$, each increasing
in $C_s$ and decreasing in $C_t$.
Then on $Q:= \{V(C_s)\cap V(C_t)=\0\}$,
\begin{equation}\label{Efg'}
\E fg \geq \E f \E g.
\end{equation}
\end{thm}

As mentioned in Section~\ref{Intro}, this can be proved along
the lines of either Theorem~\ref{1.5} or Theorem~\ref{2.5},
with Theorem~\ref{3.3} a crucial ingredient in either case.
Here we only give (sketchily) the second argument, 
leaving the reader to fill in the first (which, like the 
second, depends on Observation~\ref{U} below).

It's a little strange that we can so far prove Theorem~\ref{3.5} 
only in the absence of undirected edges, and we
conjecture that it remains true without this restriction.
(The difficulties in extending the proof below---those
for the other version are essentially the same---are 
the (related) failures of Observation~\ref{U} and 
of the validity of the hypothesis (iii) when we come to
apply Lemma~\ref{LPA}.)

\mn
{\em Proof.}
We will not repeat the proof of Theorem~\ref{2.5}, 
but just indicate what changes are needed in the present
situation.

The state space $\hOm$ and transitions for our Markov chain
are essentially as before.
(Here we have chosen to say
$C_s$, $C_t$ rather than $C_S$, $C_T$,
but as noted earlier (following Theorem~\ref{1.2})
this really makes no difference.)
Of course $\{V(C_s)\cap V(C_t)=\0\}$ now replaces
$\{S\nrar T\}$ as the conditioning event $Q$.

Let us write $\psi$ for our (unconditioned) percolation
measure, $\phi$ for our $\psi$
conditioned on $Q$ and
$\hap$ for the measure which $\phi$ induces on $\hOm$.
%We will, however, abuse this and related notation slightly:
%for an event $A$
%depending on $\om\in\{0,1\}^E$ we will sometimes use
%``$\hap$ conditioned on $A$" as shorthand for
%``the measure which $\phi$ conditioned on $A$ 
%induces on $\Om$."
The argument here then follows that for Theorem~\ref{2.5}
{\em verbatim} until, in proving positive association
of the collection (\ref{coll1}),
we come to establishing
conditions (ii) and (iii) of Lemma~\ref{LPA}.
For these we 
need the easily verified (but crucial)

\begin{obs}\label{U}
For $U\sub V$ the distribution of $C_s$ is the same under
$\phi$ conditioned on $\{V(C_t)=U\}$ as under $\psi$ 
conditioned on $\{s\nra U\}$
(and similarly with the roles of s and t reversed).
\end{obs}

In view of this, (ii) is an immediate consequence of 
Theorem~\ref{3.3}
(the relevant information from conditioning on $W$ 
being just the
resulting value of $C_s^{n-1}$).

For (iii) we first observe that if $W,W'$ 
are possible values of
$(W_1\dots W_a)$ with $W\leq W'$, and $U,U'$ are the corresponding values
of $V(C_s^{n-1})$, then according to 
Remark~\ref{inc} we have $U\sub U'$.
Moreover, by Observation~\ref{U}, 
the distribution of $(Z_1\dots Z_b)$
($=(X_e^n:e\in E)$)
given $W$ is simply the distribution 
of (the indicator of) $E\sm C_t$ under
$\psi$ conditioned on $\{t\nra U\}$.
So, writing $\hap_U$ for this distribution
on $(X_e^n:e\in E)$, 
we need to show that $U'\supseteq U$
implies $\hap_{U'}\succ \hap_{U}$
(note that increasing $C_t^n$ corresponds to 
{\em decreasing} the $X_e^n$'s).
This follows from Theorem~\ref{3.3}:
Notice that $\hap_{U'}$ is the same as $\hap_{U}$ 
conditioned on
$B:=\{t\nra U'\sm U\}$.
(More accurately, $\hap_{U'}$ is the distribution 
induced on the indicator of $E\sm C_t$ by 
$\psi(\cdot|t\nra U)$ conditioned on $B$.)
But then, since $B$ is a decreasing event determined by 
$C_t$, Theorem~\ref{3.3} says that under $\hap_{U}$, 
$B$ is negatively correlated
with any increasing event determined by $C_t$; that is, 
$\hap_{U'}\succ \hap_{U}$.\qed

\mn
{\bf Remark.}
The choice of $\hOm$ is a key to the preceding argument.  
For instance,
taking the state space to be the analogue of that in 
the alternative proof of Theorem~\ref{2.5} sketched 
at the end of Section~\ref{DMCP}---namely
$\{\om\in\{0,1\}^E:\mbox{$Q$ holds for $\om$}\}$---gets 
in trouble
because we lose some positive correlations, 
e.g. of events $\{\om_e =1\}$.

\bn
{\large\bf Acknowledgment}
We thank Tom Liggett for drawing our attention, after publication of our paper
\cite{BeKa}, to his paper \cite{BFKL} with Belitsky et al.

\end{document}